\documentstyle[12pt]{article}

\begin{titlepage}

\title{Equivariant torsion and $G$-CW-complexes}

\author{Ulrich Bunke\thanks{Mathematisches Institut, Universit\"at G\"ottingen, Bunsenstr. 3-5, 37073 G\"ottingen, GERMANY, E-mail:bunke@uni-math.gwdg.de}}

\end{titlepage}

\newcommand{\proof}{{\it Proof.$\:\:\:\:$}}

\newcommand{\R}{{\bf R}}

\newcommand{\Z}{{\bf Z}}
\newcommand{\C}{{\bf C}}

\newcommand{\res}{{\rm res}}
\newcommand{\Tr}{{\rm Tr}}

\newcommand{\cH}{{\cal H}}

\newcommand{\cE}{{\cal E}}

\newcommand{\cC}{{\cal C}}

\newcommand{\rk}{{\mbox{\rm rank}}}
\newcommand{\im}{{\mbox{\rm im}}}

\newcommand{\cF}{{\cal F}}
\newcommand{\Ree}{{\rm Re }}

\newcommand{\inter}{{\rm int}}

\newcommand{\ee}{{\rm e}}

\newcommand{\id}{{\mbox{\rm id}}}

\newcommand{\nat}{{\bf  N}}

 \newcommand{\Or}{{\rm Or }}

\def\hB{\hspace*{\fill}$\Box$ \newline\noindent}

\def\hB{\hspace*{\fill}$\Box$ \\[0.5cm]\noindent}

\newtheorem{prop}{Proposition}[section]
\newtheorem{lem}[prop]{Lemma}
\newtheorem{ddd}[prop]{Definition}

\newtheorem{kor}[prop]{Corollary}

\begin{document}

\maketitle

\tableofcontents
\newcommand{\cD}{{\cal D}}
\newcommand{\discup}{{\sqcup}}
\newcommand{\Cyl}{{\rm Cyl}}
\newcommand{\pr}{{\rm pr}}
\newcommand{\Mat}{{\rm Mat}}

\section{Introduction}

\subsection{}

In this note we consider equivariant Reidemeister
and analytic torsion invariants of closed oriented 
$G$-manifolds, where $G$ is any compact Lie group.

Equivariant analytic torsion for closed oriented
odd-dimensional $G$-manifolds 
for arbitrary compact Lie groups $G$
and equivariant Reidemeister torsion of
closed oriented $G$-manifolds for finite $G$
were introduced in \cite{lottrothenberg91} and
further studied in \cite{lueck93} and \cite{bismutzhang94}.
In the present paper we generalize the definition of equivariant
Reidemeister torsion to general compact Lie groups $G$
and address the question of equality with equivariant analytic torsion.

If $G$ is finite, then 
equivariant Reidemeister torsion is in fact an invariant
of $G$-equivariant locally constant sheaves $\cF$ of 
finite-dimensional Hilbert spaces 
over  $G$-CW-complexes. We extend this invariant
to general compact Lie groups $G$.

If $M$ is a closed $G$-manifold ($G$ a compact Lie group), then
there is a natural equivalence class of 
$G$-homotopy equivalences 
$f:X\rightarrow M$ called simple structure (see section
\ref{see22} for details), where $X$
is a $G$-CW-complex. If $\cF$ is a $G$-equivariant 
locally constant sheaf $\cF$ of 
finite-dimensional Hilbert spaces, then the equivariant 
Reidemeister torsion of $(M,\cF)$ is defined using $(X,f^*\cF)$, and it
is independent of choices.

Viewing a $G$-CW complex $X$ as a filtered $G$-space,
we express the equivariant Reidemeister torsion of $(X,\cF)$
in terms of the equivariant Reidemeister torsion of the restriction of 
$\cF$ to the $G$-cells and a contribution of 
the spectral sequence induced by the filtration.
There is a clear separation into an invariant which only depends
on restrictions of $\cF$ to the $G$-cells, and an invariant
which depends on the way the cells are glued together.
If $G$ is connected, then the latter invariant is trivial.

We compute the equivariant Reidemeister torsion and the 
equivariant analytic torsion in terms of contributions
of $G$-cells. The contribution of a $G$-cell can
be further evaluated by restricting to one-dimensional
subgroups of $G$.  
In particular we compute the equivariant Reidemeister torsion and the 
equivariant analytic torsion of compact symmetric
spaces by topological means and recover the result
of \cite{koehler97}

\subsection{}\label{r1}

Let $G$ be a compact Lie group, and $M$ be a closed odd-dimensional 
oriented $G$-manifold.
Let $F\rightarrow M$ be a $G$-equivariant flat hermitean vector bundle
and $\cF$ be the associated $G$-equivariant locally constant sheaf of 
finite-dimensional Hilbert spaces.

If we choose a $G$-equivariant Riemannian metric $g^M$, then
one defines the equivariant analytic torsion
$\rho_{an}(M,g^M,\cF):G\rightarrow \C$ (see \cite{lottrothenberg91}, \S X)
as a spectral invariant of the Laplace operator $\Delta_{g^M}$ acting on $F$-valued forms:
$$\rho_{an}(M,g^M,\cF)(g):=\frac{d}{ds}_{|s=0}\frac{1}{\Gamma(s)}\int_0^\infty (\Tr_s N g \ee^{-t\Delta_{g^M}}-\chi^\prime(M,\cF)(g)) t^{s-1 }dt\ ,$$
where $\chi^\prime(M,\cF):=\sum_{i=0}^\infty
 (-1)^i i \Tr g_{| H^i(M,\cF)}$, $N$ denotes the $\Z$-grading
of the bundle of $F$-valued forms, and the integral converges
for $\Ree(s)>>0$ and has a meromorphic continuation to all of $\C$.

By definition $\rho_{an}(M,g^M,\cF)$ is a class function on $G$.
If $\cF$ is acyclic, i.e. 
$H^*(M,\cF)=0$, then $\rho_{an}(M,g^M,\cF)$ is independent of
$g^M$, and we write $\rho_{an}(M,g^M,\cF)=:\rho_{an}(M,\cF)$.

\subsection{}\label{r2}

Let $G$ be finite, let $M$ be a closed oriented $G$-manifold,
and let $\cF$ be a $G$-equivariant locally 
constant sheaf of finite-dimensional Hilbert spaces.
We choose a $G$-Hilbert module structure on $H^*(M,\cF)$.
Using a smooth $G$-equivariant triangulation of $M$ one can define 
equivariant Reidemeister torsion
$\rho(M,\cF):G\rightarrow \C$
which is again a class function on $G$ (see \cite{lottrothenberg91}, Sec. 5).
Equivalently, there is a natural simple structure $f:X\rightarrow M$ 
and we can define $\rho(M,\cF):=\rho(X,f^*\cF)$.
We present the details of this definition in 
Section \ref{reidtors}.

Assume that $M$ is odd-dimensional.
If $\cF$ is acyclic,
then by \cite{lottrothenberg91}, Prop.16, we have the equality of class functions
\begin{equation}\label{chemuel} \rho(M,\cF)=\rho_{an}(M,\cF)\ .\end{equation}
In case that $\cF$ is not acyclic we fix a $G$-invariant 
Riemannian metric $g^M$. Let $\cH(M,g^M,\cF)$ denote the
space of harmonic $F$-valued forms. Then $\cH(M,g^M,\cF)$
is a $G$-Hilbert module. On $H^*(M,\cF)$ we choose the $G$-Hilbert module
structure such that the de Rham isomorphism
$\cH(M,g^M,\cF)\stackrel{\sim}{\rightarrow} H^*(M,\cF)$ becomes an isometry.
Again
by \cite{lottrothenberg91}, Prop.16, we have the equality
\begin{equation}\label{chemuel1} \rho(M,\cF)=\rho_{an}(M,g^M,\cF)
\end{equation} 
of class functions.
  
\subsection{}\label{r3}

Let $G$ be any compact Lie group, $M$ be a closed oriented $G$-manifold,
and $F\rightarrow M$ be a $G$-equivariant flat hermitean vector bundle.
We fix a $G$-Hilbert module structure on $H^*(M,\cF)$.

In the present subsection we define equivariant Reidemeister torsion
of $(M,\cF)$ and discuss the validity of (\ref{chemuel}) and (\ref{chemuel1}).

Let $FG:=\{g\in G\:|\: (\exists n\in\nat \:|\: g^n=1)\}$ denote the set
of elements of finite order. $FG$ is a dense subset of $G$ which is invariant
under conjugation.
Let $C(FG):=\{f:FG\rightarrow \C\:|\: f(g^h)=f(g), \forall g\in FG, h\in G\}$
denote the space of
all real-valued functions on $FG$, which are invariant under conjugation,
where $g^h:=hgh^{-1}$.

Let $\Or_f(G)$ denote the full subcategory of the orbit
category (see \cite{lueck89}, 8.16) consisting of all
homogeneous spaces $G/\Gamma$, where $\Gamma\subset G$
is finite.

We have a contravariant functor $C:\Or_f(G)\rightarrow \C-vect$ 
associating to the object $G/\Gamma$ the space of class functions 
$C(\Gamma)$ on $\Gamma$.
If $f:G/\Gamma\rightarrow  G/\Gamma^\prime$ is a morphism in 
$\Or_f(G)$, then there is a $g\in G$ such that
$\{g\gamma g^{-1}\:|\:\gamma\in\Gamma\}=\Gamma^g\subset \Gamma^\prime$ and $f(h\Gamma)=hg^{-1}\Gamma^\prime$.
The functor $C$ associates to $f$ the map
$C(f):C(\Gamma^\prime)\rightarrow C(\Gamma^g)\rightarrow C(\Gamma)$,
where the first arrow $\res^{\Gamma^\prime}_{\Gamma^g}:C(\Gamma^\prime)\rightarrow  C(\Gamma^g)$ 
is the restriction of class functions
and the second is induced by the map $\Gamma\rightarrow \Gamma^g$,
$\gamma\mapsto g\gamma g^{-1}$.

There is a natural bijection
\begin{equation}\label{ggg2}C(FG)\stackrel{\sim}{\rightarrow}\lim_{\stackrel{\Or_f(G)}{\leftarrow}} C(\Gamma)\ ,\end{equation}
which is induced by the restrictions $\res^G_\Gamma :C(FG)\rightarrow C(\Gamma)$,  $\Gamma\subset G$ finite.
 
For $\Gamma\subset G$ let $\res^G_\Gamma M$ denote the $\Gamma$-manifold
obtained from $M$ by restricting the $G$-action to $\Gamma$.
If $V$ is a $G$-module, then let $\res^G_\Gamma V$ denote the
$\Gamma$-module obtained by restriction.
Note that $H^*(\res^G_\Gamma M,\cF)=\res^G_\Gamma H^*(M,\cF)$ canonically,
and the latter has a natural $\Gamma$-Hilbert module structure.

If $G/\Gamma\in\Or_f(G)$, then $\rho(\res^G_\Gamma M,\cF)\in C(\Gamma)$
is well defined by \ref{r2}. By Proposition \ref{smore}
the collection $\{\rho(\res^G_\Gamma M,\cF)\}_{G/\Gamma\in\Or_f(G)}$
is a section of the functor $C$ and thus defines an element
$$\tilde{\rho}(M,\cF)\in\lim_{\stackrel{\Or_f(G)}{\leftarrow}} C(\Gamma)\ .$$
\begin{ddd}\label{ee}
Let $\rho(M,\cF)\in C(FG)$ be the element which corresponds to $\tilde{\rho}(M,\cF)$ under the bijection (\ref{ggg2}).
 \end{ddd}

With this definition equality of equivariant Reidemeister and analytic torsion
for arbitrary Lie groups becomes a formal consequence of the corresponding result for
finite groups. 

\begin{lem}\label{gche}
Let $G$ be a compact Lie group, $M$ be
a closed odd-dimensional oriented $G$ manifold, $g^M$ be a $G$-invariant Riemannian metric and $F\rightarrow M$
be a $G$-equivariant flat hermitean vector bundle.
We equip $H^*(M,\cF)$ with the $G$-Hilbert module structure
such that the de Rham isomorphism $\cH^*(M,g^M,\cF)\cong H^*(M,\cF)$
becomes an isometry.
Then
$$\rho(M,\cF)=\rho_{an}(M,g^M,\cF)_{|FG}\ .$$
\end{lem}
\proof
Let $g\in FG$ and $G/\Gamma\in\Or_f(G)$ be such that $g\in\Gamma$.
Then obviously 
$\rho_{an}(M,g^M,\cF)(g)=\rho_{an}(\res^G_\Gamma M,g^M,\cF)(g)$.
By (\ref{chemuel1}) we have
 $$\rho(M,\cF)(g)\stackrel{def}{=}\rho(\res^G_\Gamma M,\cF)(g)=\rho_{an}(\res^G_\Gamma M,g^M,\cF)(g)=\rho_{an}(M,g^M,\cF)(g)\ .$$
\hB

\subsection{}

Note that equivariant Reidemeister torsion depends on the choice
of a $G$-Hilbert module structure on the cohomology $H^*(M,\cF)$.
In the present subsection we show how one can define an invariant
that is independent of this choice.

Let $G^0\subset G$ denote the component of the identity of $G$.
Then we have an exact sequence
$$
 0\rightarrow G^0\rightarrow G\stackrel{q}{\rightarrow} \pi_0(G)\rightarrow 0\ .$$
Note that the restriction $q_{|FG}\rightarrow \pi_0(G)$ is still surjective.
Thus the pull-back in the sequence below defining $\hat{C}(FG)$ is injective.
$$0\rightarrow C(\pi_0(G))\stackrel{q^*}{\rightarrow} C(FG)\rightarrow \hat{C}(FG)\rightarrow 0$$
For any $G/\Gamma\in \Or_f(G)$ let $\Gamma^0:=\Gamma\cap G^0$. Then we have exact sequences
\begin{eqnarray*}
&&0\rightarrow \Gamma^0\rightarrow \Gamma\stackrel{q}{\rightarrow} \Gamma/\Gamma^0\rightarrow 0\\
&&0\rightarrow C( \Gamma/\Gamma^0)\stackrel{q^*}{\rightarrow} C(\Gamma)\rightarrow \hat{C}(\Gamma)\rightarrow 0\ ,\end{eqnarray*}
where $\hat{C}(\Gamma)$ is defined by the second sequence.

We can consider the functor
$\hat{C}:\Or_f(G)\rightarrow \C-vect$, which associates to
$G/\Gamma\in\Or_f(G)$ the space $\hat{C}(\Gamma)$.
If $f:G/\Gamma\rightarrow G/\Gamma^\prime$ is a morphism
in $\Or_f(G)$,
then the map $\hat{C}(f):\hat{C}(G/\Gamma^\prime)\rightarrow \hat{C}(G/\Gamma)$
is represented by $C(f):C(\Gamma)\rightarrow C(\Gamma^\prime)$ which maps $C(\Gamma/\Gamma^0)$ to $C(\Gamma^\prime/(\Gamma^\prime)^0)$.

Since the natural map $C(\pi_0(G))\rightarrow \lim_{\stackrel{\Or_f(G)}{\leftarrow}} C(\Gamma/\Gamma^0)$ is an
isomorphism we conclude from
$$\begin{array}{ccccccccc}
0&\rightarrow&C(\pi_0(G))&\rightarrow &C(FG)&\rightarrow&\hat{C}(FG)&\rightarrow
&0\\
&&\downarrow\|&&\downarrow\|&&\downarrow \hat{I}&&\\
0&\rightarrow&\lim_{\stackrel{\Or_f(G)}{\leftarrow}} C(\Gamma/\Gamma^0)&\rightarrow &\lim_{\stackrel{\Or_f(G)}{\leftarrow}} C(\Gamma)&\rightarrow&\lim_{\stackrel{\Or_f(G)}{\leftarrow}} \hat{C}(\Gamma)&\rightarrow& \lim^1_{\stackrel{\Or_f(G)}{\leftarrow}} C(\Gamma/\Gamma^0)
\end{array}$$
that $\hat{I}$ is injective.

Note that $G^0$ acts trivially on $H^*(M,\cF)$.
Thus if $G/\Gamma\in \Or_f(G)$, then $\Gamma^0$ acts trivially on $H^*(\res^G_\Gamma M,\cF)$, too.
By Lemma \ref{iikk} the class $\hat{\rho}(\res^G_\Gamma M,\cF)\in \hat{C}(\Gamma)$ of $\rho(\res^G_\Gamma M,\cF)\in C(\Gamma)$ is independent
of the choice of a $\Gamma$-Hilbert module structure on
$H^*(\res^G_\Gamma M,\cF)$. 

The collection $\{\hat{\rho}(\res^G_\Gamma M,\cF)\}_{G/\Gamma\in\Or_f(G)}$
defines a section of the functor $\hat{C}$ and therefore an element
$$\tilde{\hat{\rho}}(M,\cF)\in \lim_{\stackrel{\Or_f(G)}{\leftarrow}} \hat{C}(\Gamma)\ .$$
It is easy to see that $\tilde{\hat{\rho}}(M,\cF)$
is in the range of $\hat{I}$.
\begin{ddd}\label{ee1}
Let $\hat{\rho}(M,\cF)\in \hat{C}(FG)$ be the unique element such that
$\hat{I}( \hat{\rho}(M,\cF))=\tilde{\hat{\rho}}(M,\cF)$.
\end{ddd}
The class $\hat{\rho}_{an}(M,g^M,\cF)\in \hat{C}(FG)$ of 
$\rho_{an}(M,g^M,\cF)_{|FG}$ is independent of the choice of the 
$G$-invariant Riemannian
metric $g^M$ (see \cite{lottrothenberg91}, \S X).
We thus write $\hat{\rho}_{an}(M,\cF):=\hat{\rho}_{an}(M,g^M,\cF)$.
The proof of the following Lemma is similar to that of Lemma \ref{gche}.
\begin{lem}
Let $G$ be a compact Lie group, $M$ be
a closed odd-dimensional $G$ manifold, and $F\rightarrow M$
be a $G$-equivariant flat hermitean vector bundle. Then
\begin{equation}\label{gll}\hat{\rho}(M,\cF)=\hat{\rho}_{an}(M,\cF)\ .\end{equation}
\end{lem} 

\subsection{}

It is natural to ask what kind of differential-topological 
information about the $G$-manifold $M$ and the flat bundle $F$ is encoded 
in the invariants $\rho(M,\cF)$ and $\hat{\rho}(M,\cF)$.
While $\rho(M,\cF)$ contains global information about $M$ it turns
out that $\hat{\rho}(M,\cF)$ only depends on the
type of $G$-cells of $X$ and the restriction of $f^*\cF$ to the cells,
where $f:X\rightarrow M$ represents the prefered simple 
structure of the smooth closed $G$-manifold $M$.
In particular it is independent of the way the cells are glued together.

We now formulate the result in detail.
A $G$-space $G/H\times D^n$ is called a $n$-dimensional $G$-cell of
type $H$, where $H\subset G$ is a closed subgroup.
Let $X$ be a finite $G$-CW-complex and $\cF$ be a $G$-equivariant 
locally constant sheaf of finite-dimensional Hilbert spaces over $X$.  
For any $G$-cell $E=G/H_E\times D^n\hookrightarrow X$ of dimension $\dim(E):=n$ let  $\cF_E\rightarrow G/H_E$ denote the restriction of $\cF$ to $G/H_E\times\{0\}$. 
Then $\hat{\rho}(G/H_E,\cF_E)\in\hat{C}(FG)$ is defined. 
\begin{prop}[Corollary \ref{dec12}]\label{main}
Let $X$ be a finite $G$-CW-complex, and $\cF$ be a $G$-equivariant 
locally constant sheaf of finite-dimensional Hilbert spaces over $X$. 
Then we have
$$\hat{\rho}(X,\cF)=\sum (-1)^{\dim(E)}\hat{\rho}(G/H_E,\cF_E)\ ,$$
where the sum is taken over all $G$-cells of $X$.  
\end{prop}

\subsection{}

The equivariant torsion $\hat{\rho}(M,\cF)$ 
is determined by its restrictions to 
all Cartan subgroups $T$ of $G$.
If we apply Proposition \ref{main} to $\res^G_T M$, then
we can compute the equivariant torsion of $M$ in terms of the 
$T$-cells.
In Section \ref{tcell} we study the equivariant torsion of a
$T$-cell. It vanishes iff the isotropy group has codimension zero or
greater than one, and it is explicitly computable,
if the isotropy group has codimension one.

Let $f:X\rightarrow \res^G_T M$ represent the prefered simple structure.
Consider $t\in T$.
Let $I$ be the collection of  $T$-cells  $E\cong T/S_E\times D^{n_E}$ of $X$
with $\dim(T/S_E)=1$. For $E\in I$ 
let $J_E(t)$  be the
collection of connected components $E_i$, $i\in J_E(t)$, of $E$ (note that $E_i\cong S^1\times D^{n_E}$),
such that $tE_i=E_i$.
Let $m=\dim(\cF)$ and $U_i\in U(m)$ be the holonomy of
$f^*\cF_{|E_i}$. Then $U_i$ is determined uniquely by the choice of 
a base point $o_i\in E_i$, identification of the fibre
$\cF_{o_i}$ with $\C^m$, and the choice of an orientation
of $E_i$.
The matrix $U_i$ can be written as $\ee^{2\pi\imath a_i}$ for a
selfadjoint $a_i\in \Mat(m,\C)$.
The element $t$ acts as rotation of the circle-part of $E_i$ by the angle
$2\pi \tau_i$. Moreover there are unitary isomorphisms $\lambda_i$
of $\cF_{o_i}$ given by the action of $t$
composed with parallel transport back to $o_i$ 
in direction opposite to the orientation. 
Note that $\lambda_i$ and $a_i$ commute.
Then
\begin{lem}
The equivariant torsion
$\hat{\rho}(\res^G_T M,\cF)$ is given by the class of
$$T\ni t\mapsto \sum_{E\in I,i\in I_E}(-1)^{n_E}   \Tr\:\psi(\lambda_i,a_i,\tau_i)\ ,$$
where an explicit formula for $\psi$ is given in \ref{ss3}.
\end{lem}
  
For the purpose of illustration in \ref{ssym} we compute
the equivariant Reidemeister torsion of odd-dimensional symmetric
spaces. Employing (\ref{gll}) we essentially recover the results
of \cite{koehler97} about equivariant  analytic torsion
of symmetric spaces.

\section{Restriction of simple structures}

In this section we recall some basic results in equivariant topology.

\subsection{}

The $G$-space $G/H\times D^n$ is called a $n$-dimensional $G$-cell of
type $H$, where $H\subset G$ is a closed subgroup. A finite relative
$G$-CW-complex is a pair of $G$-spaces $(X,A)$ together with
a finite filtration by $G$-spaces
$$A=X_{-1}\subset X_0\subset X_1\subset \dots\subset X_N = X,\quad \cup_{n=-1}^N {X_n} = X\ ,$$
and a collection of $G$-subspaces $e^n_i\subset X_n$, $i\in I_n$, $n\ge 0$, $\sharp I_n<\infty$,
with the following properties\\
{\bf (1) :}$\:\:$  $A/G$ is Hausdorff\\
{\bf (2) :}$\:\:$ $X$ has the weak topology with respect to the filtration
$\{X_n\}$\\
{\bf (3) :}$\:\:$ for $n\ge 0$ there are $G$-push outs
$$\begin{array}{ccc}\discup_{i\in I_n} G/H_i\times S^{n-1}&\stackrel{\discup_i q_i}{\longrightarrow} &X_{n-1}\\
\downarrow&&\downarrow\\
\discup_{i\in I_n} G/H_i\times D^{n}&\stackrel{\discup_i Q_i}{\longrightarrow} &X_{n}
\end{array}\ ,
$$ such that $e^n_i=Q_i(G/H_i\times \inter \:D^n)$
(see  \cite{lueck89} Def. 1.2.).

\subsection{}\label{see22}

If $X$ is a $G$-space, then a simple structure on $X$ is given
by a pair $(Z,f)$, where $Z$ is a finite $G$-CW-complex and
$f:Z\rightarrow X$ is a $G$-homotopy equivalence.
A second pair $(Z^\prime,f^\prime)$ defines the same simple
structure on $X$, if $f^\prime_* \tau^G((f^\prime)^{-1}\circ f)=0$
holds true
in the Whitehead group $Wh^G(X)$, where $(f^\prime)^{-1}$ denotes any homotopy inverse of $f^\prime$  (see \cite{lueck89}, 4.27).

Let 
$$\begin{array}{ccc}
X_0&\rightarrow&X_1\\
\downarrow&&\downarrow\\
X_2&\rightarrow&X\end{array}$$
be a push-out of $G$-spaces and each of the $X_i$, $i=0,1,2$
be equipped with a simple structure. Then 
$X$ has a prefered simple structure \cite{lueck89}, p75.

\subsection{}

If $X$ is a closed smooth $G$-manifold (possibly with boundary), 
then it has a prefered simple
structure (\cite{lueck89}, 4.36). It is obtained by induction
over the number of orbit types.
If $X$ has one orbit type, then $X/G$ is a smooth manifold
which has a smooth triangulation. Lifting this triangulation to
$X$ we obtain a $G$-CW-decomposition of $X$ representing
the prefered simple structure.
If $X$ has several orbit types then we write $X$
as a push-out of two $G$-manifolds with less orbit types.
We apply the induction hypothesis and \ref{see22}
in order to obtain the prefered simple structure
of $X$.

\subsection{}

Let $G$ be a compact Lie group and $X$ be a finite $G$-CW complex.
If $H\subset G$ is a closed subgroup, then in general there is no
natural $H$-CW structure on $\res^G_H X$. But we have
the following result.

\begin{prop}\label{conj}
{\bf (1) :}$\:\:$ 
Let $X$ be a finite $G$-CW complex. 
If $H\subset G$ is a closed subgroup, then there exists
a prefered simple structure $f:Z\rightarrow \res^G_H X$.\\
{\bf (2) :}$\:\:$
Let $X$ be a finite $G$-CW complex. If $K\subset H\subset G$ are closed subgroups,
$f:Z\rightarrow \res^G_H X$ and $g:Y\rightarrow \res^H_K Z$ 
represent the prefered
simple structures given in {\bf (1)}, then
$f\circ g : Y\rightarrow \res^G_K X$ represents the prefered simple structure, too.\\
{\bf (3) :}$\:\:$
Let $M$ be a smooth closed $G$-manifold, and let $f:X\rightarrow M$
represent the prefered simple structure. Let $H\subset G$
be a closed subgroup, and let $g:Z\rightarrow \res^G_H X$
represent the prefered simple structure. 
Then $f\circ g : Z\rightarrow \res^G_H M$ represents the
prefered simple structure of the closed smooth $H$-manifold
$\res^G_H M$.\\
{\bf (4) :}$\:\:$
Let $X$ be a finite $G$-CW complex, $H\subset G$ be a closed subgroup,
and let $f:Z\rightarrow \res^G_H X$ represent the prefered simple structure.
Let $g\in G$, $H^g:=gHg^{-1}$, $Z^g$ be the $H^g$-CW complex
which is obtained from $Z$ by letting $u\in H^g$ act by $g^{-1}ug$,
and define $f^g:Z^g\rightarrow \res^G_{H^g} X$ by $f^g=g\circ f$.
Then $f^g:Z^g\rightarrow \res^G_{H^g} X$ represents the prefered
simple structure. \end{prop}
\proof
{\bf (1)} 
The homogeneous spaces $G/H$ for all closed
subgroups $L\subset G$ are smooth $H$-manifolds.
Thus we have prefered simple structures on the $G$-cells of $X$ considered
as $H$-spaces. Writing $X$ as a push-out
over its $G$-cells and using \ref{see22} we obtain
a prefered simple structure on $\res^G_H X$.
Note that the technical assumptions \cite{lueck89}, 7.3 and 7.23 for this procedure are satisfied (see \cite{lueck89}, 7.27, see also \cite{illman78},\cite{illman83},\cite{illman90}).\\
{\bf (2)} 
It suffices to show this for the homogeneous spaces
$G/L$. In this case we can apply {\bf (3)}
(Transitivity of the restriction was also announced in
\cite{illman90}).\\
{\bf (3)} This is \cite{lueck89}, Lemma 7.4.5.\\
{\bf (4)} It again suffices to verify this assertion for 
the homogeneous spaces $G/L$. In this case we can apply \cite{illman90}, Lemma 1.4.\hB

\newcommand{\Chi}{{\cal X}}
\newcommand{\cone}{{\rm cone}}

\section{Equivariant torsion}\label{reidtors}

\subsection{}
Let $\Gamma$ be a finite group, and let $R(\Gamma)$ denote the
representation ring of $\Gamma$ with real coefficients.
If $\pi$ is a finite-dimensional representation of $\Gamma$,
then $\chi_\pi$ denotes its character.
We have $\chi_\pi\in C(\Gamma)$, and
the map $R(\Gamma)\ni\pi\mapsto \chi_\pi\in C(\Gamma)$ induces an isomorphism of $\C$-vector spaces
$\Chi : R(\Gamma)\otimes_\R \C\cong C(\Gamma)$.
 
\subsection{}\label{alg}

Consider a finite group $\Gamma$.
Let $f:V\rightarrow W$ be an ismorphism of finite-dimensional
$\Gamma$-Hilbert modules.
If $\pi$ is an irreducible representation of $\Gamma$,
then let $V(\pi)$, $W(\pi)$ denote the $\pi$-isotypic components
and $f(\pi):V(\pi)\rightarrow W(\pi)$ the induced isomorphism.
We define $[[f]]:\in R(\Gamma)$ by
$$[[f]](\pi):= \frac{1}{2\dim(\pi)}\log |\det \: f(\pi)^*f(\pi)|\ .$$
Let 
$$\cC : \dots \rightarrow C^p\stackrel{c^p}{\rightarrow}C^{p+1}\rightarrow \dots$$
be an acyclic finite cochain complex of finite-dimensional $\Gamma$-Hilbert modules.
Then there exists
a chain contraction
$\kappa^*:C^*\rightarrow C^{*-1}$, and $c^{ev}+\kappa^{ev}:C^{ev}\rightarrow C^{odd}$ is an isomorphism, where $C^{ev}:=\oplus_{k\in\Z}C^{2k}$, $C^{odd}:=\oplus_{k\in\Z}C^{2k+1}$.
We define
$$\rho(\cC):=\Chi [[c^{ev}+\kappa^{ev}:C^{ev}\rightarrow C^{odd}]]\ .$$
Note that $\rho(\cC)$ does not depend on the choice of the chain
contraction $\kappa$.

Let $f:\cC\rightarrow \cD$ be a chain homotopy equivalence.
Then we consider the complex $\cone(f)$ with
$\cone(f)^n:=C^n\oplus D^{n-1}$ and the differential
$$\left( \begin{array}{cc} c^n & 0\\ f^n & - d^{n-1}\end{array}\right)\ .$$
We define $t(f)=\rho(\cone(f))\in C(\Gamma)$.
Note that if $g:\cC\rightarrow \cD$ is homotopy equivalent to $f$, then
$t(f)=t(g)$.
If $g:\cD\rightarrow \cE$ is a second chain homotopy equivalence,
then $t(g\circ f)=t(f)-t(g)$.

If $\cC$ is a finite complex of finite-dimensional $\Gamma$-Hilbert modules,
then we consider the complex $\cH(\C)$ with $p$'th space
$H^p(\cC)$ and trivial differential.
If $f:\cC\rightarrow\cD$ is a homotopy equivalence, then we obtain a 
homotopy equivalence $f_*:\cH(\cC)\rightarrow \cH(\cD)$. 

Let $\cC$, $\cD$ be finite complexes of finite-dimensional $\Gamma$-Hilbert modules with prefered $\Gamma$-Hilbert module structures on $H^*(\cC)$, $H^*(\cD)$. A homotopy equivalence $f:\cC\rightarrow
\cD$ is called simple if $t(f)+t(f_*)=0$.

We call two finite complexes of finite-dimensional $\Gamma$-Hilbert modules $\cC,\cD$ with prefered $\Gamma$-Hilbert module structures on $H^*(\cC)$, $H^*(\cD)$
equivalent, if there exists a simple homotopy equivalence $f:\cC\rightarrow
\cD$. We write $[\cC]$ for the equivalence class.

If $\cC$ is a finite complex of finite-dimensional $\Gamma$-Hilbert modules
with prefered $\Gamma$-Hilbert module structure on $H^*(\cC)$,
and if $i:\cH(\cC)\rightarrow \cC$ is any embedding such that $i_*=\id$,
then we set $\rho(\cC):=-t(i)$.
If $\cC$ and $\cD$ are equivalent, then $\rho(\cC)=\rho(\cD)$,
hence we can write
$\rho([\cC]):=\rho(\cC)$, where $\cC$ is any representative of
$[\cC]$.

Let $\Gamma^0\subset \Gamma$ be a subgroup such that $\Gamma^0$
acts trivially on $\cH(\cC)$. Let $\hat{C}(\Gamma):=C(\Gamma)/C(\Gamma/\Gamma^0)$
and $\hat{\rho}(\cC)$ 
be the class of $\rho(\cC)$ in $\hat{C}(\Gamma)$.
\begin{lem}\label{iikk}
$\hat{\rho}(\cC)$ is independent of the choice of the 
$\Gamma$-Hilbert module structure on $H^*(\cC)$.
\end{lem}
\proof
Let $\cH_j(\cC)$, $j=1,2$, be the complex $\cH(\cC)$ equipped
with two $\Gamma$-Hilbert module structures.
Let $\rho_j(\cC)$ be the corresponding torsion.
Then we have
$\rho_2(\cC)=\rho_1(\cC)-t(\id:\cH_2(\cC)\rightarrow \cH_1(\cC))$.
It is easy to see that $t(\id:\cH_2(\cC)\rightarrow \cH_1(\cC))\in C(\Gamma/\Gamma^0)$. \hB

\subsection{}\label{ttff}
Let $\Gamma$ be a finite group and
$X$ be a $\Gamma$-space of the homotopy type of a finite $\Gamma$-CW-complex.
Let $\cF$ be a $\Gamma$-equivariant locally constant sheaf of finite-dimensional
Hilbert spaces. We fix a $\Gamma$-Hilbert module structure on $H^*(X,\cF)$.

Let $(Z,f)$ represent a simple structure on $X$.
Then we form the cellular cochain complex
$\cC(Z,f^*\cF)$, which is a finite complex of finite-dimensional $\Gamma$-Hilbert modules. We equip $H^*(\cC(Z,f^*\cF))$ with the
$\Gamma$-Hilbert module structure such that the canonical map $f^*:H^*(X,\cF)\rightarrow H^*(\cC(Z,f^*\cF))$ becomes an isometry.

If $(Z^\prime,f^\prime)$ and $(Z,f)$ represent the same simple structure
of $X$,
then $\cC(Z,f^*\cF)$ and $\cC(Z^\prime,(f^\prime)^*\cF)$
are equivalent chain complexes.
If $X$ is a $\Gamma$-space with distinguished simple structure and $\Gamma$-Hilbert module structure on $H^*(X,\cF)$,
then we write $[\cC(X,\cF)]:=[\cC(Z,f^*\cF)]$,
where $(Z,f^*\cF)$ is any representative of the distinguished simple
structure of $X$. We define
$$\rho(X,\cF):=\rho([\cC(X,\cF)])\in C(\Gamma)\ .$$

Let $\Gamma^0\subset \Gamma $ act trivially on $H^*(X,\cF)$. Then
the class $\hat{\rho}(X,\cF)\in\hat{C}(\Gamma)$ does not depend
on the choice of the $\Gamma$-Hilbert module structure
on $H^*(X,\cF)$.

\subsection{}\label{s1}

Let $\Gamma$ be a finite group and
$\Gamma^\prime\subset \Gamma$.
 If $Z$ is a $\Gamma$-CW complex, then since $\Gamma$ is finite $\res^\Gamma_{\Gamma^\prime} Z$
carries a natural $\Gamma^\prime$-CW structure.
Moreover $\id:\res^\Gamma_{\Gamma^\prime} Z\rightarrow \res^\Gamma_{\Gamma^\prime} Z$ represents  the prefered simple structure
given in Proposition \ref{conj} {\bf (1)}.

Let $X$ be a $\Gamma$-space of the homotopy type of a finite $\Gamma$-CW complex.
If $(Z,f)$ represents a simple structure for $X$,
then $(\res^\Gamma_{\Gamma^\prime} Z,f)$ represents a simple structure
of $\res^\Gamma_{\Gamma^\prime} X$.

We choose a $\Gamma$-Hilbert module structure on $H^*(X,\cF)$
which we also use for $H^*(\res^\Gamma_{\Gamma^\prime} X,\cF)=\res^\Gamma_{\Gamma^\prime}H^*(X,\cF)$.

\begin{lem}\label{resres}
$$\res^\Gamma_{\Gamma^\prime} \rho(X,\cF)= \rho( \res^\Gamma_{\Gamma^\prime}X,\cF)\ .$$
\end{lem}
\proof
We use $(\res^\Gamma_{\Gamma^\prime} Z,f)$ to represent the prefered simple structure of $\res^\Gamma_{\Gamma^\prime} X$.
Then
\begin{eqnarray}\cC(\res^\Gamma_{\Gamma^\prime} Z,f^*\cF)&=&\res^\Gamma_{\Gamma^\prime} \cC(Z,f^*\cF)\nonumber\\
{}[\cC(\res^\Gamma_{\Gamma^\prime} Z,f^*\cF)]&=&[\res^\Gamma_{\Gamma^\prime} \cC(Z,f^*\cF)]\label{mnb}\ .
\end{eqnarray}
Let $h:V\rightarrow W$ be an isomorphism of finite-dimensional Hilbert-$\Gamma$-modules
and $\res^\Gamma_{\Gamma^\prime} h:\res^\Gamma_{\Gamma^\prime} V
\rightarrow \res^\Gamma_{\Gamma^\prime} W$.
Then for any irreducible representation $\tau$ of $\Gamma^\prime$ we have
\begin{eqnarray*}
[[\res^\Gamma_{\Gamma^\prime} h]](\tau)&=&\frac{1}{2\dim(\tau)}
\log|\det (\res^\Gamma_{\Gamma^\prime} h(\tau)^* \res^\Gamma_{\Gamma^\prime} h(\tau)|\\
&=&\sum_{\pi\in\hat{\Gamma}} \frac{[\pi:\tau]}{2\dim(\pi)}\log|\det h^*(\pi)h(\pi)|\\
&=&\sum_{\pi\in\hat{\Gamma}} [\pi:\tau] [[h]](\pi)
\end{eqnarray*}
Note that $\res^\Gamma_{\Gamma^\prime}\chi_\pi=\sum_{\tau\in\hat{\Gamma^\prime}} [\pi:\tau] \chi_\tau$.
Thus
\begin{eqnarray}
\Chi([[\res^\Gamma_{\Gamma^\prime} h]])&=&\sum_{\tau\in\hat{\Gamma^\prime}} [[\res^\Gamma_{\Gamma^\prime} h]](\tau)\chi_\tau\label{mnb1}\\
&=&\sum_{\tau\in\hat{\Gamma^\prime}} \sum_{\pi\in\hat{\Gamma}} [\pi:\tau] [[h]](\pi)  \chi_\tau\nonumber\\
&=&\sum_{\pi\in\hat{\Gamma}}   [[h]](\pi)   \res^\Gamma_{\Gamma^\prime}\chi_\pi\nonumber\\
&=&\res^\Gamma_{\Gamma^\prime} \Chi([[h]])\nonumber\ .
\end{eqnarray}
 The Lemma now follows from (\ref{mnb}) and (\ref{mnb1}). 
\hB 

\subsection{}\label{s2}

Let $G$ be a compact Lie group,
$X$ be a finite $G$-CW-complex
and $\cF$ be a $G$-equivariant locally constant sheaf of finite-dimensional
Hilbert spaces. 
We choose a $G$-Hilbert module structure on $H^*(X,\cF)$
which induces a $\Gamma$-Hilbert module structure on $H^*(\res^G_\Gamma X,\cF)=\res^G_\Gamma H^*(X,\cF)$ for all $G/\Gamma\in \Or_f(G)$.

If $G/\Gamma\in\Or_f(G)$, then $\res_\Gamma^G X$
has a prefered simple structure by Proposition \ref{conj} {\bf (1)} and
$\rho(\res_\Gamma^G X,\cF)\in C(\Gamma)$ is defined.
If $h\in G$, then let 
$\Gamma^h:=h\Gamma h^{-1}$, $G/\Gamma^h\in\Or_f(G)$, and 
set $g^h:=hgh^{-1}$ for  $g\in G$. 
\begin{lem}\label{conj1}
If $g\in\Gamma$ and $h\in G$, then $\rho(\res_\Gamma^G X,\cF)(g)=\rho(\res_{\Gamma^h}^G X,\cF)(g^h)$.
\end{lem}
\proof
Let $(Z,f)$ represent the prefered simple structure of $\res^G_\Gamma X$.
Then by Proposition \ref{conj} {\bf (4)} (we employ the
notation introduced there) the pair $(Z^h,f^h)$
represents the prefered simple structure of $\res^G_{\Gamma^h} X$.
We have an isomorphism of complexes of $\Gamma^h$-Hilbert modules
$C(Z,f^*\cF)^h=C(Z^h,(f^h)^*\cF)$, where
the underlying space of $C(Z,f^*\cF)^h$ is $C(Z,f^*\cF)$
and $\Gamma^h$ acts by $g^h\mapsto g$.  
Similarly we have isomorphisms of complexes $\Gamma^h$-Hilbert modules
$\cH(\res^G_\Gamma X,\cF)^h=\cH(\res^G_{\Gamma^h} X,\cF)$.
If $i,i^h,j$ denote the inclusions 
$i:\cH(\res^G_\Gamma X,\cF)\hookrightarrow C(Z,f^*\cF)$,
$i^h:\cH(\res^G_\Gamma X,\cF)^h\hookrightarrow C(Z,f^*\cF)^h$,
$j:\cH(\res^G_{\Gamma^h} X,\cF)\hookrightarrow C(Z^h,(f^h)^*\cF)$,
then we have $t(i)(g)=t(i^h)(g^h)=t(j)(g^h)$.
This proves the Lemma. \hB

\subsection{}

We keep the assumptions of \ref{s2}
Then the Lemmas \ref{conj1} and \ref{resres} imply 

\begin{kor}\label{rrrr}
The collection $\{\rho(\res^G_\Gamma X,\cF)\}_{G/\Gamma\in\Or_f(G)}$,
defines a section of the functor $C$ (see Subsection \ref{r3}).
\end{kor}

Let now $M$ be a closed $G$-manifold and $F\rightarrow M$
be a $G$-equivariant flat hermitean vector bundle.
We fix $G$-Hilbert module structures on $H^*(M,\cF)$.
If $G/\Gamma\in \Or_f(G)$, then $\res^G_\Gamma M$ has a prefered simple structure. We employ this structure in order to define $\rho(\res^G_\Gamma M,\cF)$
as explained in Section \ref{ttff}.

\begin{prop}\label{smore}
The collection $\{\rho(\res^G_\Gamma M,\cF)\}_{G/\Gamma\in\Or_f(G)}$,
defines a section of the functor $C$.
\end{prop}
\proof
Let $f:X\rightarrow M$ represent the prefered simple structure.
Then by Proposition \ref{conj}, {\bf (3)}, 
we have $\rho(\res^G_\Gamma M,\cF)=\rho(\res^G_\Gamma X,f^*\cF)$.
Thus the Proposition is implied by Corollary \ref{rrrr}. \hB

\section{Computations}

\subsection{}\label{fil}

Let $G$ be a compact Lie group and $X$ be a finite $G$-CW-complex.
Let $\cF$ be a $G$-equivariant locally constant sheaf of finite-dimensional
Hilbert spaces. We choose a $G$-Hilbert module structure on $H^*(X,\cF)$.
Being a $G$-CW complex $X$ has a natural filtration
$\emptyset=X_{-1}\subset X_0\subset X_1\subset \dots\subset X_N=X$.
Consider $G/\Gamma\in  \Or_f(G)$.
\begin{lem}
There exists a representative $f:Z\rightarrow \res^G_\Gamma X$
of the prefered simple structure such that\\
{\bf (1) :}$\:\:$
$Z$ is filtered by $\Gamma$-CW subcomplexes
$\emptyset=Z_{-1}\subset Z_0\subset Z_1\subset \dots\subset Z_N=Z$
and $f_{|Z_p}:Z_p\rightarrow \res^G_\Gamma X_p$ represents the prefered
simple structure for all $p\in\{0,1,\dots,N\}$.\\
{\bf (2)}
If 
$$\dots\subset F_{p+1}\cC(Z,f^*\cF)\subset F_{p}\cC(Z,f^*\cF)\subset\dots$$
denotes the decreasing filtration of the associated cochain complexes
(we write $F_p:=F_p\cC(Z,f^*\cF)$), then
$F_p/F_{p+1}$ is the cochain complex associated to
a representative of the prefered simple structure of $\discup_{i\in I_p} \res^G_\Gamma G/H_i\times (D^p,S^{p-1})$ and the local system $Q^*\cF$,
where $Q$ is given by $$Q:=\discup_{i\in I_p}Q_i:\discup_{i\in I_p}  G/H_i\times D^p \rightarrow X_p$$
(recall that $I_p$ is the indexing set of the $p$-dimensional $G$-cells of $X$
and $Q_i$ denote chacteristic maps). 
\end{lem}
\proof
The construction of $f:Z\rightarrow \res^G_\Gamma X$ goes by induction and is based on \cite{lueck89}, 4.29-4.32.
For $X_{-1}$ the assertion is trivial.
Assume that we have constructed 
the simple structure $(Z_{n-1},f_{n-1})$ of $\res^G_\Gamma X_{n-1}$
together with the filtration by $\Gamma$-subspaces $(Z_{n-1})_m$.
Then we have to construct a simple structure $(Z_n,f_n)$ of $\res^G_\Gamma X_n$ together with the filtration by $\Gamma$-subspaces $(Z_n)_m$.

Let $X_n$ be given by the $G$-push out
$$\begin{array}{ccc}\discup_{i\in I_n} G/H_i\times S^{n-1}&\stackrel{\discup_i q_i}{\longrightarrow} &X_{n-1}\\
\downarrow&&\downarrow\\
\discup_{i\in I_n} G/H_i\times D^{n}&\stackrel{\discup_i Q_i}{\longrightarrow} &X_{n}
\end{array}\ .
$$ 
We choose representatives $h_i:V_i\rightarrow \res^G_\Gamma G/H_i$ of the prefered simple structures. Then $h_i\times\id:U_i:=V_i\times D^n\rightarrow \res^G_\Gamma G/H_i\times D^n$
is a simple structure on the cell $\res^G_\Gamma G/H_i\times D^n$.
We choose a $\Gamma$-equivariant cellular map $p:\discup_{i\in I_n} U_i\times S^{n-1}\rightarrow Z_{n-1}$ such that $f_{n-1}\circ p\sim_\Gamma \discup_{i\in I_n} q_i \circ l_i$
and $\sim_\Gamma$ stands for $\Gamma$-homotopic. 
We now replace $(Z_{n-1},f_{n-1})$ by $(\Cyl(p),f^\prime_{n-1})$,
which represents the same simple structure on $X_{n-1}$ ($f^\prime_{n-1}$ has still to be constructed).
The filtration of $\Cyl(p)=\discup_{i\in I_n} U_i\times S^{n-1}\times I\cup_p Z_{n-1}$ is given by
$\Cyl(p)_{m}=(Z_{n-1})_m\subset \Cyl(p)$ for $m\le n-2$.
Let $\pr:\Cyl(p)\rightarrow Z_{n-1}$ be the projection.
In order to construct $f_{n-1}^\prime$ we consider the following $\Gamma$-homotopy commutative diagram $$\begin{array}{ccc}
\discup_{i\in I_n}\res^G_\Gamma G/H_i\times S^{n-1}\discup X_{n-2}&\stackrel{\discup_i  q_i\discup j}{\rightarrow}& \res^G_\Gamma X_{n-1}\\
\uparrow\discup_i l_i\discup f_{n-2}&&\uparrow f_{n-1}\circ \pr\\
\discup_{i\in I_n} U_i\times S^{n-1}\discup \Cyl(p)_{n-2}&\stackrel{J}{\hookrightarrow}&\Cyl(p)\end{array}\ ,$$
where $j:X_{n-2}\hookrightarrow X_{n-1}$ is the inclusion,
$l_i:=(h\times \id)_{|U_i\times S^{n-1}}$, and 
$J_{|\discup_{i\in I_n} U_i\times S^{n-1}}$ is the natural identification with
the closed subspace 
$\discup_{i\in I_n} U_i\times S^{n-1}\times \{0\}\subset\discup_{i\in I_n} U_i\times S^{n-1}\times I$. 
Since $J$ is a cofibration we can find $f_{n-1}^\prime\sim_\Gamma  f_{n-1}\circ \pr$ such that 
$f_{n-2}\circ J_{|\Cyl(p)_{n-2}}=(f^\prime_{n-1})_{|\Cyl(p)_{n-2}}$ and the following diagram commutes:
$$\begin{array}{ccccc}
\discup_{i\in I_n} \res^G_\Gamma G/H_i\times D^{n}&\hookleftarrow &\discup_{i\in I_n} \res^G_\Gamma G/H_i\times S^{n-1} &\stackrel{\discup_i q_i}{\rightarrow}&\res^G_\Gamma X_{n-1}\\
\uparrow &&\uparrow&&\uparrow f^\prime_{n-1} \\
\discup_{i\in I_n} U_i\times D^{n}&\hookleftarrow&\discup_{i\in I_n} U_i\times S^{n-1}&\stackrel{J}{\hookrightarrow}&\Cyl(p)\end{array}\ .$$
Let $Z_{n}$ be the $\Gamma$-push out
$$\begin{array}{ccc}
\discup_{i\in I_n} U_i\times S^{n-1}& \rightarrow & \Cyl(p)\\
\downarrow&&\downarrow\\
\discup_{i\in I_n} U_i\times D^n&\rightarrow &Z_n  \end{array}$$
and $f_n:Z_n\rightarrow \res^G_\Gamma X_n$ be the natural map of push outs.
Then $f_n:Z_n\rightarrow \res^G_\Gamma X_n$ represents the prefered simple structure,
and $Z_n$ is filtered by $\Gamma$-subspaces $(Z_n)_m$ such that $(f_n)_{|(Z_n)_m}:(Z_n)_m\rightarrow \res^G_\Gamma X_m$
represents the prefered simple structure for all $m\le n$.
This finishes the proof of {\bf (1)}. Assertion {\bf (2)} is an easy consequence of the construction. \hB
 
\subsection{}

For $i\in I_p$ define $$\bar{Q}_i:G/H_i\cong G/H_i\times\{0\} \hookrightarrow G/H_i\times D^p \stackrel{Q_i}{\rightarrow} X$$
and set $\cF_i:=\bar{Q}_i^*\cF$.
For any $p\in\nat_0$ and $i\in I_p$ we fix $G$-Hilbert module
structures on $H^*(G/H_i, \cF_i)$.
This induces $G$-Hilbert module structures on 
the cohomology complexes $\cH(G/H_i, \cF_i)$.
If $\cC$ is a cochain complex, then let $\cC[p]$ be the cochain complex
with $\cC[p]^n:=\cC^{n+p}$ obtained from $\cC$.
We equip $H^*(F_p/F_{p+1})$
with $\Gamma$-Hilbert module
structures such that the natural isomorphism $\cH(F_p/F_{p+1})\cong
\oplus_{i\in I_p}\cH(G/H_i,\cF_i)[p]$ becomes an isometry.
Then we have
$$[F_p/F_{p-1}] =\oplus_{i\in I_p} [\cC(U_i,h_i^*\cF_i)[p]]\ .$$

We are going to express $\rho(\res^G_\Gamma X,\cF)$ in terms of
$\rho(\res^G_\Gamma G/H_i,\cF_i)$ and a contribution of the spectral sequence
$\cE:=(E^{p,q}_r,d_r)$ associated to the filtration of $X$.
For trivial $\Gamma$ this was worked out
in \cite{lueckschickthielmann96}. But
\cite{lueckschickthielmann96}, Lemmas 4.6 and 4.7, extend immediately
to the case of a finite group $\Gamma$. We recall the
result. Let
\begin{eqnarray*}
Z_r^{p,q}&:=&\im(H^{p+q}(F_p/F_{p+r})\rightarrow H^{p+q}(F_p/F_{p+1}))\\
B_r^{p,q}&:=&\im(H^{p+q-1}(F_{p-r+1}/F_{p})\rightarrow H^{p+q}(F_p/F_{p+1}))\\
E_r^{p,q}&:=&Z_r^{p,q}/B_r^{p,q}\\
Z_\infty^{p,q}&:=&\im(H^{p+q}(F_p)\rightarrow H^{p+q}(F_p/F_{p+1}))\\
B_\infty^{p,q}&:=&\im(H^{p+q-1}(\cC(Z,f^*\cF))\rightarrow H^{p+q}(F_p/F_{p+1}))\\
E_\infty^{p,q}&:=&Z_\infty^{p,q}/E_\infty^{p,q}\\
F^{p,q}&:=&\im(H^{p+q}(F_p)\rightarrow H^{p+q}(\cC(Z,f^*\cF)))
\ .
\end{eqnarray*}
There are natural isomorphisms $\psi^{p,q}:F^{p,q}/F^{p+1,q-1}\rightarrow E^{p,q}_\infty$.

Note that $H^*(\cC(Z,f^*\cF))\stackrel{f^*}{\cong} H^*(X,\cF)$
has a prefered $\Gamma$-Hilbert module structure.
We equip $Z_r^{p,q}, B_r^{p,q}, E_r^{p,q}$, $Z_\infty^{p,q}$, $B_\infty^{p,q}$, $E_\infty^{p,q}$, $F^{p,q}$
with the corresponding (sub)quotient $\Gamma$-Hilbert module
structures.
For any $p,q,r$ we have a complex of finite $\Gamma$-Hilbert modules
$$\cE_r^{p,q}: \dots\rightarrow E_r^{p+nr,q-(r-1)n}\rightarrow E_r^{p+(n+1)r,q-(r-1)(n+1)}\rightarrow \dots\ ,$$
and $H^n(\cE_r^{p,q})=\cE_{r+1}^{p+nr,q-(r-1)n}$ has a prefered $\Gamma$-Hilbert module structure.
Thus $\rho(\cE_r^{p,q})$ is well defined.
Furthermore note that
$\rho(F_p/F_{p+1})=(-1)^p\sum_{i\in I_p} \rho(\res^G_\Gamma G/H_i,\cF_i)$.
The following Proposition can be proved by repeating the argument of the proof
of \cite{lueckschickthielmann96}, Thm. 4.4.
\begin{prop}\label{dec}
$$\rho(\res^G_\Gamma X,\cF)=\sum_{p} (-1)^p\sum_{i\in I_p} \rho(\res^G_\Gamma G/H_i,\cF_i)+\sum_{r\ge 1}\sum_{p=0}^{r-1}\sum_{q} (-1)^{p+q}\rho(\cE_r^{p,q})-\sum_{p,q}(-1)^{p+q}\Chi[[\psi^{p,q}]]\ .$$
\end{prop}
Note that all three terms of the right-hand side may 
depend on the choice of the prefered $\Gamma$-Hilbert 
module structures on $H^*(G/H_i,\cF_i)$, while the left hand side
does not.

\subsection{}

Let $X$ be a finite $G$-CW-complex and $\cF$ a $G$-equivariant 
locally constant sheaf of finite-dimensional Hilbert spaces. 
Let $G/\Gamma \in\Or_f(G)$. Then $\hat{\rho}(\res^G_\Gamma X,\cF)\in \hat{C}(\Gamma)$ is well defined. 

\begin{prop}\label{dec1}
$$\hat{\rho}(\res^G_\Gamma X,\cF)=\sum_{p} (-1)^p \sum_{i\in I_p}\hat{\rho}(\res^G_\Gamma G/H_i,\cF_i)$$
\end{prop}
\proof
We fix $G$-Hilbert module structures on $H^*(X,\cF)$
and $H^*(G/H_i,\cF_i)$.
Let $q:\Gamma\rightarrow \Gamma/\Gamma^0$ be the projection.
We have to show that 
\begin{equation}\label{zuei}\rho(\cE_r^{p,q})\in q^* C(\Gamma/\Gamma^0), \quad \Chi[[\psi^{p,q}]]\in q^* C(\Gamma/\Gamma^0)\ .\end{equation}

Viewing $X$ as a filtered $G$-space we see that there is a spectral sequence
$\tilde{\cE}:=(\tilde{E}_r^{p,q},d_r)$ of $G$-modules with
\begin{eqnarray*}
\tilde{Z}_r^{p,q}&:=&\im(H^{p+q}(X_{p+r-1},X_{p-1})\rightarrow \oplus_{i\in I_{p}}H^{p+q}(G/H_i,\cF_i))\ ,\\
\tilde{B}_r^{p,q}&:=&\im(H^{p+q-1}(X_{p-1}/X_{p-r})\rightarrow \oplus_{i\in I_{p}}H^{p+q}(G/H_i,\cF_i))\ ,\\
\tilde{E}_r^{p,q}&:=&\tilde{Z}_r^{p,q}/\tilde{B}_r^{p,q}\\
\tilde{Z}_\infty^{p,q}&:=&\im(H^{p+q}(X,X_{p-1},\cF)\rightarrow  \oplus_{i\in I_{p}}H^{p+q}(G/H_i,\cF_i))\\
\tilde{B}_\infty^{p,q}&:=&\im(H^{p+q-1}(X,\cF))\rightarrow  \oplus_{i\in I_{p}}H^{p+q}(G/H_i,\cF_i))\\
\tilde{E}_\infty^{p,q}&:=&\tilde{Z}_\infty^{p,q}/\tilde{E}_\infty^{p,q}\\
\tilde{F}^{p,q}&:=&\im(H^{p+q}(X,X_{p-1},\cF)\rightarrow H^{p+q}(X,\cF))
\end{eqnarray*}
and $G$-equivariant maps $\tilde{\psi}^{p,q}:\tilde{F}^{p,q}/\tilde{F}^{p+1,q-1}\rightarrow \tilde{E}_\infty^{p,q}$
such that $\cE=\res^G_\Gamma \tilde{\cE}$, $\psi^{p,q}=\res^G_\Gamma \tilde{\psi}^{p,q}$.
Consider the exact sequence
$$0\rightarrow G^0\rightarrow G\stackrel{q}{\rightarrow} \pi_0(G)\rightarrow 0\  .$$
Observe that the representation of $G$ on $\tilde{\cE}$ and $\tilde{F}^{p,q}$  factors over $q$. In particular $\Gamma^0$ is represented trivially.
Thus (\ref{zuei}) follows. \hB

Fix $G$-Hilbert module structures on $H^*(X,\cF)$
and $H^*(G/H_i,\cF_i)$. Then Proposition \ref{dec1} has the
\begin{kor}\label{dec12}
$$\hat{\rho}(X,\cF)=\sum_{p} (-1)^p \sum_{i\in I_p}\hat{\rho}(G/H_i,\cF_i)\ .$$
\end{kor}

\section{Reduction to Cartan subgroups}
\newcommand{\Aut}{{\rm Aut}}
\newcommand{\bC}{{\bf D}}

\subsection{}\label{ss1}

Let $G$ be a compact Lie group, $X$ be a finite $G$-CW complex, and let $F\rightarrow X$ be an equivariant flat hermitean vector bundle.
 
Recall that a Cartan subgroup $T$ of $G$ is a topologically cyclic
closed subgroup such that the Weyl group $N_G(T)/T$ is finite,
where $N_G(T)$ is the normalizer of $T$ in $G$.
If $T\subset G$ is a Cartan subgroup then it is isomorphic to the
product of a torus and a finite cyclic group
(see \cite{broeckertomdieck85}, p.177 ff, for more
details about Cartan subgroups).

 If $g\in G$, then there exists a Cartan subgroup containing $g$. The conjugacy classes of Cartan subgroups
are in natural bijection with the cyclic subgroups of $\pi_0(G)$.

Let $\{T_i\}_{i\in C}$ be a set of representatives of conjugacy classes
of Cartan subgroups of $G$.
Then restriction defines inclusions
$$\oplus_{i\in C} \res^G_{T_i}:C(FG)\rightarrow \oplus_{i\in C} C(FT_i),\quad \oplus_{i\in C} \hat{\res}^G_{T_i}:\hat{C}(FG)\rightarrow \oplus_{i\in C} \hat{C}(FT_i)\ .$$
In order to determine $f\in \hat{C}(FG)$ it is thus sufficient to 
compute $\hat{\res}^G_{T_i}f \in \hat{C}(FT_i) $ for all $i\in C$.

Let $T$ be any Cartan subgroup of $G$ and let $f:Z\rightarrow \res^G_T X$
represent the prefered simple structure. Then $\hat{\res}^G_{T}\hat{\rho}(X,\cF) = \hat{\rho}(Z,f^*\cF)$.

\subsection{}\label{tcell}

We further study the contribution of the $T$-cells of $Z$.
Let $S\subset T$ be any closed subgroup, 
and let $\cF$ be a $T$-equivariant local
system on $T/S$. 
Then $\rho(T/S,\cF)\in \hat{C}(FT)$ is well defined.
\begin{lem}
If $\dim(T/S)\not= 1$, then $\hat{\rho}(T/S,\cF)\in \hat{C}(FT)=0$.
\end{lem}
\proof
Assume first that $T/S$ is even-dimensional. Then $T/S$ is orientable.
Let $t\in FT$ generate the finite group $H\subset \Aut(T/S)$.
Then $H$ acts by orientation-preserving diffeomorphisms on $T/S$.
We employ \cite{lueck93}, Prop. 3.23, which says that
$\hat{\rho}(\res^T_H T/S,\cF)$ can be derived from 
the Poincar\'e torsion $\rho_{pd}^H(T/S,\cF)$ of the $H$-manifold
$\res^T_H T/S$ (see loc.cit. Definition 3.19).
Since $T$ is abelian, we see that $H$ acts freely or trivially on $T/S$.
In both cases $\rho_{pd}^H(T/S,\cF)$ vanishes (see loc.cit. Prop. 3.20), and thus $\hat{\res}^T_H\hat{\rho}(T/S,\cF)=0$.
Since $t$ was arbitrary we conclude that $\hat{\rho}(T/S,\cF)=0$.

Now we consider the case that $T/S$ is odd-dimensional.
We fix a $T$-invariant Riemannian metric $g^{T/S}$ on $T/S$
and equip $H^*(T/S,\cF)$ with the $T$-module structure such that the de Rham isomorphism becomes an isometry.
Then $\rho(T/S,\cF)\in C(FT)$ is well-defined.

By Lemma \ref{gche} we have $\rho(T/S,\cF)=\rho_{an}(T/S,g^{T/S},\cF)_{|FT}$.
If $\dim(T/S)>1$, then we find two everywhere linearly
independend unit-length Killing vector fields in $C^\infty(T/S,T(T/S))$.
Using the induced decomposition
of the de Rham complex by a standard argument $\rho_{an}(T/S,g^{T/S},\cF)=0$.
This implies the lemma.
\hB

\subsection{}\label{kl12}

Let $S\subset T$ be any closed subgroup such that $\dim(T/S)=1$,
and let $\cF$ be a $T$-equivariant locally constant sheaf of Hilbert spaces.
Since $T$ is abelian and $\rho(T/S,\cF)$ is additive with respect to $\cF$
without loss of generality we can assume 
that $\cF$ is one-dimensional.

We fix an orientation on $T/S$. 
The space $T/S$ is a disjoint union
of oriented circles $C_1,\dots,C_s$. Let $\ee^{2\pi\imath a_i}\in S^1$ denote
the holonomy of $\cF$ on $C_i$. We choose the Riemannian metric $g^{T/S}$ such that
the circles have length $1$.

Let $t\in T$. 
We distinguish two cases\\ 
a) that $tC_i=C_i$ for all $i$\\
b) that $tC_i\not=C_i$ for all $i$.\\
In case b) we have $\rho_{an}(T/S,g^{T/S},\cF)(t)=0$.
In case a) we have
$\rho_{an}(T/S,g^{T/S},\cF)(t)=\sum_{i=1}^s\rho_{an}(C_i,g^{C_i},\cF_i)(t)$,
where $\cF_i$ is the restriction of $\cF$ to $C_i$.

Thus we are reduced to the following situation.
Let $C_i$ be a circle, $F\rightarrow C_i$
be a one-dimensional flat hermitean bundle over $C_i$, and $t$ be
an automorphism of $F\rightarrow C_i$, which acts as rotation on $C_i$.
We parametrize $C_i=\R/\Z$ with coordinate $x$. Then $t(x)=x+\tau_i$, $\tau_i\in \R/\Z$.
Sections of $F_{|C_i}$ are identified with functions $f:\R\rightarrow \C$
satisfying $f(x+1)=\ee^{2\pi\imath a_i} f(x)$.
Then $(tf)(x)=\lambda_i f(x-\tau_i)$ for certain $\lambda_i\in S^1$.

We have $\rho(C_i,g^{C_i},\cF)_{an}(t)=\psi(\lambda_i,a_i,\tau_i)$,
and $\psi$ will be determined in Subsection \ref{ss3}.
We obtain $\rho_{an}(T/S,g^{T/S},\cF)(t)=\sum_{i=1}^s \psi(\lambda_i,a_i,\tau_i)$.
If $\cF$ is higher-dimensional, then $\lambda_i$ and $a_i$ become
commuting diagonalizable matrices, and we have 
$\rho_{an}(T/S,g^{T/S},\cF)(t)=\sum_{i=1}^s \Tr\:\psi(\lambda_i,a_i,\tau_i)$.

\subsection{}\label{ss3}

In this Subsection we derive a formula for $\psi(\lambda,a,\tau)$. 
Let $C$ be a circle, $F\rightarrow C$
be a one-dimensional flat hermitean bundle over $C$, and $t$ be
an automorphism of $F\rightarrow C$, which acts as rotation on $C$.

We parametrize $C=\R/\Z$ with coordinate $x$. Then $t(x)=x+\tau$, $\tau\in \R/\Z$.
Sections of $F$ are identified with functions $f:\R\rightarrow \C$
satisfying $f(x+1)=\ee^{2\pi\imath a} f(x)$.
Then $(tf)(x)=\lambda f(x-\tau)$ for certain $\lambda\in S^1$.

We identify one-forms on $C$ with functions using the basis
$dx\in C^\infty(C,T^*C)$.
Let $\Delta_F=-(d/dx)^2$ be the Laplace operator.
The eigenvectors of $\Delta_F$ are $f_n(x)=\exp(2\pi\imath (n+a)  x)$,
and the corresponding eigenvalue is $\mu_n=4\pi^2 (n+a)^2$.
The action of $t$ on the eigenspace spanned  by $f_n$ is
multiplication by $\lambda \ee^{- 2\pi\imath (n +a)\tau}$.

First assume that $a\in (0,1)$.
Let 
$$F_n(s):=\frac{\lambda \ee^{-2\pi\imath (n +a)\tau}}{\Gamma(s)}\int_0^\infty  t^{s-1} \ee^{-4\pi^2 (n+ a)^2 t} dt=\frac{\lambda \ee^{-2\pi\imath (n +a)\tau}}{4^s\pi^{2s} (n+ a)^{2s}}\ ,$$
then
$$\psi(\lambda,a,\tau)=-\frac{d}{ds}_{|s=0} \sum_{n\in\Z}F_n(s)\ ,$$
where we employ a meromorphic continuation of the sum $\sum_{n\in\Z}F_n(s)$
which converges for $\Ree(s)>1$.

Let $\phi(y,a,s):=\sum_{n=0}^\infty \frac{\ee^{2\pi\imath y n}}{(n+a)^s}$.  Then
$$\psi(\lambda,a,\tau)=-\frac{d}{ds}_{|s=0}\frac{\lambda}{ 4^s\pi^{2s} }\left( \ee^{-2\pi\imath \tau a} 
\phi(-\tau,a,2s)+ \ee^{-2\pi\imath \tau (a-1)} \phi(\tau,1-a,2s)\right)\ .$$

If $a=0$ and $\tau\in (0,1)$, then with $\phi(y,s):=\sum_{n=1}^\infty \frac{\ee^{2\pi\imath y n}}{n^s}$
we obtain
$$\psi(\lambda,0,\tau)=-\frac{d}{ds}_{|s=0}\frac{\lambda}{4^s\pi^{2s}} (\phi(\tau,2s)+\phi(-\tau,2s))\ .$$
Using \cite{lottrothenberg91}, Prop. 31, one can show that
$$\psi(\lambda,0,\tau)=\lambda \bC + \lambda\left(\frac{\Gamma(\tau)^\prime}{\Gamma(\tau)}+\frac{\Gamma(1-\tau)^\prime}{\Gamma(1-\tau)}\right)$$
for some  explicitly known constant $\bC\in \C$.
In the case that $a=0$ and $\tau=0$ we have
$$\psi(\lambda,0,0)=-\frac{d}{ds}_{|s=0} \frac{2\lambda\zeta_R(2s)}{4^s\pi^{2s}}\ ,$$
where $\zeta_R$ denotes the Riemann zeta function.

\subsection{}

In this Subsection we combine the results of Subsections \ref{ss1}, \ref{tcell}, and \ref{kl12}.
For simplicity we assume that $G$ is a connected compact Lie group. 
Let $X$ be a finite $G$-CW complex, and let $\cF\rightarrow X$ be a $G$-equivariant locally constant sheaf of finite-dimensional Hilbert spaces. 
Let $T\subset G$ be a maximal torus, and 
let $f:Z\rightarrow \res^G_T X$ represent the prefered simple structure.
 
Let $I$ be the index set for the $T$-cells $E=T/S_E\times D^{n_E}$ of $Z$ 
with $\dim(T/S_E)=1$. 
Since $T$ is connected the quotient $T/S$ is a circle $S^1$.
If $t\in T$, then for each $E\in I$ we fix an orientation
of $T/S_E$, and we define the
rotation number $\tau_E(t)\in \R/\Z$, the constant $\lambda_E(t)$, and the holonomy $\ee^{2\pi\imath a_E}$ of $\cF_{E}$.
We obtain the following Proposition.

\begin{prop}\label{lll3}
$\hat{\rho}(X,\cF)$ is uniquely determined by its restriction to $T$,
which is represented by the function
$$FT\ni t\mapsto \sum_{E\in I} (-1)^{n_E}\Tr \: \psi(\lambda_E(t),a_E,\tau_E(t))\ .$$
\end{prop}

\subsection{}\label{ssym}

In this section we compute the equivariant Reidemeister
torsion of odd-dimensional symmetric spaces. We recover
results of \cite{koehler97} using a completely
different method.
Let $\theta$ be the constant sheaf with fibre
$\C$, and let $\hat{\rho}(M):=\hat{\rho}(M,\theta)$
for any closed $G$-manifold $M$, where we equip
$\theta$ with the obvious $G$-action.

Let $G/K$ be a compact, irreducible, odd-dimensional symmetric
space. We assume that $G$ is connected.
\begin{lem}\label{lll1}
If $\hat{\rho}(G/K)\not=0$, then $\rk G=\rk K+1$,
and 
$$G/K=\left\{\begin{array}{c}SO(2m)/SO(2p-1)\times SO(2m-2p+1)\:\:\:\mbox{or}\\
                              SU(3)/SO(3)\end{array} \right. \ .$$
\end{lem}
\proof
Let $T$ be a maximal torus of $G$
and $f:X\rightarrow \res^G_T (G/K)$ be a representative of the
prefered simple structure. Using the construction given in
\cite{lueck89}, 4.36, we can choose $X$ such that
it has the same set of $T$-orbit types as $G/K$.
If $\hat{\rho}(G/K)=\hat{\rho}(X)\not=0$, then 
by Proposition \ref{lll3}
there exists a one-dimensional $T$-orbit $TgK\subset G/K$.
Then $T^{g^-1}\cap K$ is a $\rk G-1$-dimensional torus in $K$.
Hence $\rk G\ge \rk K\ge \rk G -1$. If $\rk K= \rk G$, then
$G/K$ is even-dimensional. Thus $\rk K = \rk G -1$.
The second assertion follows from the classification of irreducible 
compact symmetric spaces. \hB

Now assume that $G/K$ is a compact, irreducible, odd-dimensional symmetric
space with $\rk K = \rk G - 1$.

\begin{lem}\label{lll2}
There is a one-to-one correspondence of one-dimensional
$T$-orbits in $G/K$ and $W_G(T)/W_K(T)$
given by $ N_G(T)\ni g\mapsto Tg^{-1}K$,
where $W_G(T):=N_G(T)/T$ and $W_K(T):=N_K(T)/T\cap K$.
\end{lem}
\proof
We can assume that $T\cap K=S$ is a maximal torus of $K$.
If $T^\prime \subset G$ is a second maximal torus of $G$
with $T\cap K=S$, then $T^\prime=T$.
Indeed, on the level of Lie algebras
we have $t^\prime=s\oplus (k^\perp)^S=t$.

Let $TgK$ be a one-dimensional $T$-orbit in $G/K$.
Then $T^{g^-1}\cap K$ is a maximal torus in $K$, and hence 
$T^{g^-1}\cap K=S^k$ for a suitable $k\in K$.
Replacing $g$ by $gk^{-1}$ we obtain
$T^{g^-1}\cap K=S$, and thus $g^{-1}\in N_G(T)$.
If $g\in T\cup K$, then $TgK=TK$.
Thus the correspondence 
of $W_G(T)/W_K(T)$ with the set of one-dimensional
orbits is established by $g^{-1}\in N_G(T)\mapsto TgK$. \hB

Let $s\subset t$ be the Lie algebra of $S\subset T$.
Let $a:=t/s$, and let $L\subset a$
be the lattice of those $[l]\in t/s$, $l\in t$, which satisfy
$\exp(l)\in S$. We identify $a\cong \R$ such that
$L$ is identified with $\Z$. 
The exponential map yields an identification
$i:\R/\Z \cong a/L\cong T/S$. 
If $t\in T$, then let $\alpha(t)\in \R/\Z$ be such that
$t i(x) = i(\alpha(t)+x)$, $\forall x\in \R$.

For $\tau\in\R/\Z$, $\tau\not=0$, define 
$$\psi(\tau):=\bC+\left(\frac{\Gamma(\hat{\tau})^\prime}{\Gamma(\hat{\tau})}+\frac{\Gamma(1-\hat{\tau})^\prime}{\Gamma(1-\hat{\tau})}\right)\ ,$$
where $\hat{\tau}\in (0,1)$ represents $\tau$.
If $\tau =0$, we put $$\psi(\tau):=-\frac{d}{ds}_{|s=0} \frac{2\lambda\zeta_R(2s)}{4^s\pi^{2s}}\ .$$

Because of the equality (\ref{gll}) the following proposition recovers the computation of equivariant analytic torsion
of compact symmetric spaces \cite{koehler97}, Thm. 11, up to
a constant function.

\begin{prop}
Let $G/K$ be a compact, irreducible, odd-dimensional symmetric
space.
If $\rk K\not= \rk G - 1$, then $\hat{\rho}(G/K)=0$.
If $\rk K= \rk G - 1$, then let
$T\subset G$ be a maximal torus such that $S:=T\cap K$ is a maximal
torus of $K$. Then the restriction of
$\hat{\rho}(G/K)$ to $T$ is represented by
$$FT\ni t\mapsto \frac{1}{\sharp W_K(T)}\sum_{w\in W_G(T)} \psi(\alpha(t^w))\ .$$
\end{prop}
\proof
Let $f:X\rightarrow \res^G_{N_G(T)} (G/K)$ be a representative of the
natural simple structure. By Lemma \ref{lll2} there is exactly one isolated one-dimensional
$N_G(T)$-orbit in $G/K$. Constructing $X$ by the inductive procedure 
given in \cite{lueck89}, 4.36, we can assume that
$X$ has exactly one cell $E=N_G(T)/N_K(T)\times D^{n_E}$ with
one-dimensional $N_G(T)$-orbits. 
Moreover $n_E=0$. 

The $T$-space $\res^{N_G(T)}_T (X)$
has a natural $T$-CW structure (since $N_G(T)/T$ is finite),
and $\res^{N_G(T)}_T (E)$ is the only $T$-cell with one-dimensional $T$-orbits.
But $\res^{N_G(T)}_T (E)$ is the disjoint union of spaces $T/S^{g}$, $g\in W_G(T)/W_K(T)$.
The proposition now follows from Proposition \ref{lll3}. 
\hB\\

{\bf Acknowledgement:}{\em The author thanks W. L\"uck
for his interest in this work and many valuable hints.}

\bibliographystyle{plain}

\end{document}